\documentclass[11pt,reqno]{amsart}
\usepackage{amsmath,amsthm,amssymb,verbatim,enumerate,ifthen,times}
\usepackage[mathscr]{eucal}
\usepackage[utf8]{inputenc}
\newcommand{\N}{\mathbb{N}}
\newcommand{\Z}{\mathbb{Z}}
\newcommand{\R}{\mathbb{R}}

\newcommand{\CM}{\mathscr{CM}}
\newcommand{\CP}{\mathscr{CP}}
\newcommand{\A}{\mathscr{A}}
\newcommand{\h}{\mathscr{H}}
\newcommand{\G}{\mathscr{G}}
\newcommand{\B}{\mathscr{B}}

\def\d{\,{\rm d}}

\newtheorem{theorem}{Theorem}
\newtheorem*{theorem*}{Theorem}
\def\Thm#1#2{\ifthenelse{\equal{#1}{*}}{\begin{theorem*}#2\end{theorem*}}
  {\begin{theorem}\label{T#1}#2\end{theorem}}}
\newtheorem{Atheorem}{Theorem}

\def\thm#1{Theorem~\ref{T#1}}

\newtheorem{proposition}[theorem]{Proposition}
\newtheorem*{proposition*}{Proposition}
\def\Prp#1#2{\ifthenelse{\equal{#1}{*}}{\begin{proposition*}#2\end{proposition*}}
             {\begin{proposition}\label{P#1}#2\end{proposition}}}

\newtheorem{corollary}[theorem]{Corollary}
\newtheorem*{corollary*}{Corollary}
\def\Cor#1#2{\ifthenelse{\equal{#1}{*}}{\begin{corollary*}#2\end{corollary*}}
             {\begin{corollary}\label{C#1}#2\end{corollary}}}

\newtheorem{lemma}[theorem]{Lemma}
\newtheorem*{lemma*}{Lemma}
\def\Lem#1#2{\ifthenelse{\equal{#1}{*}}{\begin{lemma*}#2\end{lemma*}}
             {\begin{lemma}\label{L#1}#2\end{lemma}}}
\def\lem#1{Lemma~\ref{L#1}}

\newtheorem{example}[theorem]{Example}
\newtheorem*{example*}{Example}
\def\Exa#1#2{\ifthenelse{\equal{#1}{*}}{\begin{example*}\rm #2\end{example*}}
             {\begin{example}\label{Ex#1}\rm #2\end{example}}}

\newtheorem{problem}[theorem]{Problem}

\theoremstyle{definition}
\newtheorem{definition}[theorem]{Definition}

\newtheorem{remark}[theorem]{Remark}
\newtheorem*{remark*}{Remark}
\def\Rem#1#2{\ifthenelse{\equal{#1}{*}}{\begin{remark*}\rm #2\end{remark*}}
             {\begin{remark}\label{R#1}\rm #2\end{remark}}}

\newcommand{\eq}[1]{\eqref{E#1}}
\newcommand{\Eq}[2]{\ifthenelse{\equal{#1}{*}}
  {\begin{equation*}\begin{aligned}[]#2\end{aligned}\end{equation*}}
  {\begin{equation}\begin{aligned}[]\label{E#1}#2\end{aligned}\end{equation}}}
\newcommand{\DET}[1]{\begin{vmatrix}#1\end{vmatrix}}

\long\def\comment#1{}
\def\Det#1#2{\left|\!\begin{array}{cc}#1&#2\end{array}\!\right|}

\begin{document}
\large

\date{\today}

\title[Invariance of nonsymmetric Bajraktarevi\'c mean]
{On the invariance equation for two-variable weighted nonsymmetric Bajraktarevi\'c means}

\author[Zs. P\'ales]{Zsolt P\'ales}
\address{Institute of Mathematics, University of Debrecen,
H-4002 Debrecen, Pf.\ 400, Hungary}
\email{pales@science.unideb.hu}

\author[A. Zakaria]{Amr Zakaria}
\address{Department of Mathematics, Faculty of Education, Ain Shams University, Cairo 11341, Egypt}
\email{amr.zakaria@edu.asu.edu.eg}


\dedicatory{Dedicated to the 70th birthday of Professor Karol Baron}

\thanks{The research of the first author was supported by the Hungarian Scientific Research Fund (OTKA) Grant
K-111651 and by the EFOP-3.6.1-16-2016-00022 project. This project is co-financed by the European Union
and the European Social Fund.}

\begin{abstract}
The purpose of this paper is to investigate the invariance of the arithmetic mean with respect to two
weighted Bajraktarevi\'c means, i.e., to solve the functional equation
\Eq{*}{
  \bigg(\frac{f}{g}\bigg)^{\!\!-1}\!\!\bigg(\frac{tf(x)+sf(y)}{tg(x)+sg(y)}\bigg)
  +\bigg(\frac{h}{k}\bigg)^{\!\!-1}\!\!\bigg(\frac{sh(x)+th(y)}{sk(x)+tk(y)}\bigg)=x+y \qquad(x,y\in I),
}
where $f,g,h,k:I\to\R$ are unknown continuous functions such that $g,k$ are nowhere zero on $I$, the ratio functions $f/g$,
$h/k$ are strictly monotone on $I$, and $t,s\in\R_+$ are constants different from each other. By the main result of this paper, the solutions of the above invariance equation can be expressed either in terms of hyperbolic functions or in terms of trigonometric functions and an additional weight function. For the necessity part of this result, we will assume that $f,g,h,k:I\to\R$ are four times continuously differentiable.
\end{abstract}
\maketitle

\section{Introduction}

Throughout this paper, the symbols  $\R$, and $\R_+$ will stand for the sets of
real, and positive real numbers, respectively, and $I$ will always denote a nonempty open real interval. The
classes of continuous strictly monotone and continuous positive real-valued functions defined on $I$ will be
denoted by $\CM(I)$ and $\CP(I)$, respectively.

In the sequel, a function $M:I^2\to I$ is called a
\textit{two-variable mean} on $I$ if the following so-called mean value property
\Eq{1}{
  \min(x,y)\leq M(x,y)\leq \max(x,y)   \qquad(x,y\in I)
}
holds. Also, if both of the inequalities in \eq{1} are strict whenever $x\neq y$, then we say that $M$ is a \textit{strict mean} on $I$. The \textit{arithmetic} and \textit{geometric} means are well known instances for strict means on $\R_+$. More generally, if $p$ is a
real number, then the \textit{two-variable H\"older mean} $H_p:\R_+^2\to\R$ is defined as
\Eq{*}{
  H_p(x,y)
     :=\begin{cases}
       \left(\dfrac{x^p+y^p}{2}\right)^{\frac{1}{p}}
          &\mbox{ if } p\neq0,\\
       \sqrt{xy}&\mbox{ if } p=0
       \end{cases}\qquad (x,y\in\R_+).
}
A classical generalization of H\"older means is the notion of \textit{two-variable quasi-arithmetic mean} (cf.\
\cite{HarLitPol34}), which is introduced as follows: For a continuous strictly monotone function $f:I\to\R$, the \textit{two-variable quasi-arithmetic mean}
$A_f:I^2\to I$ is defined by
\Eq{*}{
   A_f(x,y):=f^{-1}\left(\frac{f(x)+f(y)}{2}\right)\qquad (x,y\in I).
}
Given two continuous strictly monotone functions $\varphi_1,\varphi_2:I\to\R$, the \textit{generalized quasi-arithmetic means} $A_{(\varphi_1,\varphi_2)}:I^2\to I$ is defined by
\Eq{*}{
   A_{(\varphi_1,\varphi_2)}(x,y):=(\varphi_1+\varphi_2)^{-1}\left(\frac{\varphi_1(x)+\varphi_2(y)}{2}\right)\qquad (x,y\in I),
}
which was introduced by Matkowski \cite{Mat10b}.
For parameters $p,q\in\R$, the \textit{two-variable Gini mean} $G_{p,q}:\R_+^2\to\R_+$ is defined by
\Eq{*}{
   G_{p,q}(x,y)
      :=\begin{cases}
         \bigg(\dfrac{x^p+y^p}{x^q+y^q}\bigg)^{\frac1{p-q}}&\mbox{if }p\neq q,\\[4mm]
         \exp\bigg(\dfrac{x^p\log(x)+y^p\log(y)}{x^p+y^p}\bigg)&\mbox{if }p=q,
        \end{cases}
    \qquad(x,y\in\R_+).
}
(See \cite{Gin38}.) The \textit{two-variable Stolarsky mean} $S_{p,q}:\R_+^2\to\R_+$ is defined for $x,y\in\R_+$, by
\Eq{*}{
   S_{p,q}(x,y)
      :=\begin{cases}
         \bigg(\dfrac{q(x^p-y^p)}{p(x^q-y^q)}\bigg)^{\frac1{p-q}}&\mbox{if }pq(p-q)(x-y)\neq0,\\[4mm]
         \exp\bigg(-\dfrac1p+\dfrac{x^p\log(x)-y^p\log(y)}{x^p-y^p}\bigg)&\mbox{if }p=q, pq(x-y)\neq0,\\[4mm]
         \bigg(\dfrac{x^p-y^p}{p(\log(x)-\log(y))}\bigg)^{\frac1{p}}&\mbox{if }q=0, p(x-y)\neq0,\\[4mm]
         \bigg(\dfrac{x^q-y^q}{q(\log(x)-\log(y))}\bigg)^{\frac1{q}}&\mbox{if }p=0, q(x-y)\neq0,\\[4mm]
         \sqrt{xy}                                                  &\mbox{if } p=q=0,\\[4mm]
         x                                                          &\mbox{if }x=y.
        \end{cases}
}
(See \cite{Sto75}.)

Given three strict means $M,N,K:\R_+^2\to\R_+$, we say that the triple $(M,N,K)$ satisfies the \textit{invariance equation} if
\Eq{2}{
K(M(x,y),N(x,y))=K(x,y)  \qquad(x,y\in \R_+)
}
holds. If \eq{2} is valid, then we say that $K$ is invariant with respect to the mean-type mapping $(M,N)$.
The easiest example when the invariance equation is satisfied is the well-known identity
\Eq{*}{
\sqrt{xy}=\sqrt{\frac{x+y}{2}\cdot\frac{2xy}{x+y}}   \qquad(x,y\in \R_+).
}
The last identity means that
\Eq{*}{
\G(x,y)=\G(\A(x,y),\h(x,y))     \qquad(x,y\in\R_+),
}
where $\A$, $\G$, and $\h$ are the two-variable arithmetic, geometric, and harmonic means, respectively. Another invariance equation is the identity
\Eq{*}{
\A\otimes\G(x,y)=\A\otimes\G(\A(x,y),\G(x,y))       \qquad(x,y\in\R_+),
}
where $\A\otimes\G$ denotes Gauss's arithmetic-geometric mean and defined as follows
\Eq{*}{
\A\otimes\G(x,y)=\bigg(\frac{2}{\pi}\int_0^\frac{\pi}{2}\frac{\d t}{\sqrt{x^2\cos^2t+y^2\sin^2t}}\bigg)^{-1}     \qquad(x,y\in \R_+).
}

The invariance equation in more general classes of means has been studied by several authors in a large number of papers.
The invariance equation of H\"older mean solved completely by Dar\'oczy and P\'ales \cite{DarPal02c}.

\Thm{*}{
Let $p,q,r\in\R$. Then the invariance equation
\Eq{*}{
  H_r\big(H_p(x,y),H_q(x,y)\big)=H_r(x,y)
    \qquad(x,y\in\R_+)
}
is satisfied if and only if one of the following two possibilities holds:
\begin{enumerate}[(i)]
\item $p=q=r$, i.e., all all the three means are equal to each other,
\item $p+q=r=0$, i.e., $H_r$ is the geometric mean and $H_p=H_{-q}$.
\end{enumerate}
}

The more general invariance equation for quasi-arithmetic means was first solved under infinitely many times differentiability by Sut\^o \cite{Sut14a}, \cite{Sut14b} and later by Matkowski \cite{Mat99a} under twice continuous differentiability. Without imposing unnecessary regularity conditions, this problem was finally solved by Dar\'oczy and P\'ales \cite{DarPal02c}.

\Thm{*}{
Let $f,g,h:I\to\R$ be continuous strictly monotone functions. Then the invariance equation
\Eq{*}{
A_f\bigl( A_g(x,y),A_h(x,y)\bigr)=A_f(x,y)\qquad(x,y\in I)
}
holds if and only if there exist $a,b,c,d,p\in\R$ with $ac\neq0$ such that
\Eq{*}{
  g=aE_p\circ f+b,\qquad h=cE_{-p}\circ f+d,
}
where
\Eq{*}{
 E_p(t):=\begin{cases}
              \exp{(pt)}&\mbox{if } p\neq0,\\ t&\mbox{if } p=0.
           \end{cases}
}}

Burai \cite{Bur07a} and Jarczyk-Matkowski \cite{JarMat06} studied the invariance equation involving three weighted arithmetic means.  Jarczyk \cite{Jar07} solved this problem without additional regularity assumptions. The invariance of the arithmetic mean with respect to Lagrangian mean has been investigated by Matkowski (cf.\ \cite{Mat05}). The invariance of the arithmetic, geometric, and harmonic means has been studied by Matkowski \cite{Mat04c}. The following result of Baj\'ak and P\'ales \cite{BajPal09a} describes the invariance of the arithmetic mean with respect to generalized quasi-arithmetic means.

\Thm{*}{
Let $\varphi_1,\varphi_2,\psi_1,\psi_2:I\to\R$ be four times continuously differentiable functions such that $\varphi'_1\varphi'_2$ is positive on $I$. Then the functional equation
\Eq{*}{
(\varphi_1+\varphi_2)^{-1}\left(\frac{\varphi_1(x)+\varphi_2(y)}{2}\right)
 +(\psi_1+\psi_2)^{-1}\left(\frac{\psi_1(x)+\psi_2(y)}{2}\right)=x+y\qquad (x,y\in I),
}
 holds if and only if
\begin{enumerate}[(i)]
\item either there exist real constants $p,a_1,a_2,c_1,c_2,b_1,b_2,d_1,d_2$ with $p\neq0$, $a_{1}a_{2}>0$, $c_{1}c_2>0$ and $a_{1}c_1=a_{2}c_2$ such that for $x\in I$,
\Eq{*}{
\varphi_1(x)&=a_{1}E_p(x)+b_1,\qquad &\varphi_2(x)&=a_{2}E_p(x)+b_2,\\
\psi_1(x)&=c_{1}E_{-p}(x)+d_1,\qquad &\psi_2(x)&=c_{2}E_{-p}(x)+d_2;
}
\item or there exist real constants $a,b,c,d_1,d_2$ with $ac\neq0$ such that, for $x\in I$,
\Eq{*}{
\varphi_1(x)+\varphi_2(x)=ax+b,\qquad \psi_1(x)=c\varphi_2(x)+d_1,\qquad\mbox{and}\qquad\psi_2(x)=c\varphi_1(x)+d_2.
}
\end{enumerate}
}

Recently, Baj\'ak and P\'ales \cite{BajPal09b}, \cite{BajPal10} have solved the invariance equations of two-variable Gini and Stolarsky means, respectively.

\Thm{*}{
Let $p,q,r,s,u,v\in\R$. Then the invariance equation
\Eq{*}{
G_{p,q}(G_{r,s}(x,y),G_{u,v}(x,y))=G_{p,q}(x,y)\qquad (x,y\in\R_+)
}
holds if and only if one of the following possibilities holds
\begin{enumerate}[(i)]
\item $p+q=r+s=u+v=0$, that is, all the three means are equal to the geometric mean,
\item $\{p,q\}=\{r,s\}=\{u,v\}$, that is, all the three means are equal to each other,
\item $\{r,s\}=\{-u,-v\}$ and $p+q=0$, that is, $G_{p,q}$ is the geometric mean and $G_{r,s}=G_{-u,-v}$,
\item there exist $a,b\in\R$ such that $\{r,s\}=\{a+b,b\}$,
$\{u,v\}=\{a-b,-b\}$, and $\{p,q\}=\{a,0\}$, in this case $G_{p,q}$ is the power mean,
\item there exists $c\in\R$ such that $\{r,s\}=\{3c,c\}$, $u+v=0$, and $\{p,q\}=\{2c,0\}$, in this case $G_{p,q}$ is a power mean and $G_{u,v}$ is the geometric mean,
\item there exists $c\in\R$ such that $r+s=0$, $\{u,v\}=\{3c,c\}$, and $\{p,q\}=\{2c,0\}$, in this case $G_{p,q}$ is a power mean and $G_{r,s}$ is the geometric mean.
\end{enumerate}
}
\Thm{*}{
Let $p,q,r,s,u,v\in\R$. Then the invariance equation
\Eq{*}{
S_{u,v}(S_{p,q}(x,y),S_{r,s}(x,y))=S_{u,v}(x,y)\qquad (x,y\in\R_+)
}
holds if and only if one of the following possibilities holds
\begin{enumerate}[(i)]
\item $p+q=r+s=u+v=0$, that is, all the three means are equal to the geometric mean,
\item $\{p,q\}=\{r,s\}=\{u,v\}$, that is, all the three means are equal to each other,
\item $\{p,q\}=\{-r,-s\}$ and $u+v=0$, that is, $G_{u,v}$ is the geometric mean and $S_{p,q}=S_{-r,-s}$.
\end{enumerate}
}

Now we recall the notion of weighted two-variable Bajraktarevi\'c mean, the class of means where we are going to solve the invariance problem of the arithmetic mean.

Given two continuous functions $f,g:I\to\R$ such that $g$ is nowhere zero on $I$ and
the ratio function $f/g$ is strictly monotone on $I$, the \textit{weighted two-variable Bajraktarevi\'c mean}
$B_{f,g}:I^2\times\R_+^2\to I$ is defined by
\Eq{*}{
  B_{f,g}(x,y;t,s):=\bigg(\frac{f}{g}\bigg)^{\!\!-1}\!\!\bigg(\frac{tf(x)+sf(y)}{tg(x)+sg(y)}\bigg)
  \qquad(x,y\in I;\,t,s\in\R_+).
}
See the paper \cite{Baj58} for the original definition. These means have been extensively investigated by Acz\'el--Dar\'oczy \cite{AczDar63c}, Dar\'oczy--Losonczi \cite{DarLos70}, Losonczi \cite{Los71a}, \cite{Los71b}, \cite{Los71c}, \cite{Los99}, Losonczi--P\'ales \cite{LosPal08}, \cite{LosPal11a}, P\'ales \cite{Pal87d}. For recent generalizations of Bajraktarevi\'c means, the solutions of the comparison, equality and homogeneity problems, we refer to the papers of the authors \cite{PalZak17} and \cite{PalZak18a}.

The purpose of this paper is to investigate the invariance of the arithmetic mean with respect to two
weighted Bajraktarevi\'c means, i.e., to solve
\Eq{inv}{
  B_{f,g}(x,y;t,s)+B_{h,k}(x,y;s,t)=x+y \qquad(x,y\in I)
}
or equivalently,
\Eq{*}{
  \bigg(\frac{f}{g}\bigg)^{\!\!-1}\!\!\bigg(\frac{tf(x)+sf(y)}{tg(x)+sg(y)}\bigg)
  +\bigg(\frac{h}{k}\bigg)^{\!\!-1}\!\!\bigg(\frac{sh(x)+th(y)}{sk(x)+tk(y)}\bigg)=x+y,
}
where $f,g,h,k:I\to\R$ are continuous functions such that $g,k$ are nowhere zero on $I$, the ratio functions $f/g$ and
$h/k$ are strictly monotone on $I$, and $t,s\in\R_+$ are constants different from each other.

For the sake of convenience and brevity, we introduce certain regularity classes as follows. Let the class
$\B_0(I)$ contain all pairs $(f,g)$ such that
\begin{enumerate}
 \item[(a)] $f,g:I\to\R$ are continuous functions,
 \item[(b)] $g$ is nowhere zero on $I$, and
 \item[(c)] $f/g$ strictly monotone on $I$.
\end{enumerate}
For $n\geq1$, let $\B_n(I)$ denote the class of all pairs $(f,g)$ such that
\begin{enumerate}
 \item[(+a)] $f,g:I\to\R$ are $n$ times continuously differentiable functions,
 \item[(b)] $g$ is nowhere zero on $I$, and
 \item[(+c)] $(f/g)'$ is nowhere zero on $I$.
\end{enumerate}
Obviously, condition (+a) and (+c) imply (a) and (c), respectively. Therefore, $\B_n(I)\subseteq\B_0(I)$ for
all $n\geq1$ and thus the Bajraktarevi\'c mean $B_{f,g}$ is well-defined for $(f,g)\in\B_n(I)$. The following
lemma explains the need of the regularity classes $\B_n(I)$.

\Lem{reg}{
Let $n\in\N$ and $(f,g)\in\B_n(I)$. Then, for every $t,s\in\R_+$, the mapping
\Eq{*}{
  (x,y)\mapsto B_{f,g}(x,y;t,s) \qquad(x,y\in I)
}
is $n$-times continuously differentiable on $I^2$.
}

The proof is an easy consequence of standard calculus rules and therefore is left for the reader.
In what follows, we establish an implicit equation for the definition of weighted two-variable Bajraktarevi\'c
means.

\Lem{sol}{
Let $(f,g)\in\B_0(I)$. Then, for every $x,y\in I$ and $t,s\in\R_+$, the value
$z:=B_{f,g}(x,y;t,s)$ is the unique solution of the equation
\Eq{d}{
  \left|\begin{array}{cc}
  tf(x)+sf(y) & f(z) \\ tg(x)+sg(y) & g(z)
  \end{array}\right|=0.
}
}
The proof of this lemma is again elementary, hence it is omitted.

We say that two pairs of functions $(f,g)$ and $(h,k)$ are equivalent (and we write $(f,g)\sim(h,k)$)
if there exist constants $a,b,c,d$ with $ad\neq cd$ such that
\Eq{*}{
  h=af+bg \qquad\mbox{and}\qquad k=cf+dg.
}

The equivalence is a key property for characterizing the equality of weighted two-variable Bajraktarevi\'c
means.

\Lem{equ0}{
Let $(f,g),(h,k)\in\B_0(I)$. Then the equality
\Eq{*}{
  B_{f,g}(x,y;t,s)=B_{h,k}(x,y;t,s) \qquad(x,y\in I,\,t,s\in\R_+)
}
holds if and only if $(f,g)\sim(h,k)$.
}

To obtain another important characterization of the equivalence of pairs $(f,g),(h,k)\in\B_n(I)$ when
$n\geq2$, for $i,j\in\{1,\dots,n\}$ and $(f,g)\in\B_n(I)$, we define
\Eq{DD}{
  W^{i,j}_{f,g}:=\left|\begin{array}{cc} f^{(i)} & f^{(j)} \\ g^{(i)} & g^{(j)} \end{array}\right|
  \qquad\mbox{and}\qquad
  \Phi_{f,g}:=\frac{W^{2,0}_{f,g}}{W^{1,0}_{f,g}},\qquad
  \Psi_{f,g}:=-\frac{W^{2,1}_{f,g}}{W^{1,0}_{f,g}}.
}
Also, it should be noted that, for all $i,j\in\{1,\dots,n\}$,
\Eq{W1}{
W^{i,j}_{f,g}=-W^{j,i}_{f,g},\qquad W^{i,i}_{f,g}=0.
}

\Lem{equ2}{
Let $(f,g),(h,k)\in\B_2(I)$. Then $(f,g)\sim(h,k)$ holds if and only if
\Eq{*}{
  \Phi_{f,g}=\Phi_{h,k}\qquad\mbox{and}\qquad
  \Psi_{f,g}=\Psi_{h,k}.
}}

The following lemma clarifies the importance the auxiliary functions $\Phi_{f,g}$ and $\Psi_{f,g}$.

\Lem{DE}{Let $(f,g)\in\B_2(I)$. Then $f,g$ are solutions of the second-order differential equation
\Eq{DE}{
  y''=\Phi_{f,g}y'+\Psi_{f,g}y.
}}

\begin{proof}
Using the definitions of $\Phi_{f,g},\Psi_{f,g}$ from \eq{DD}, we can rewrite equation \eq{DE} in the following equivalent
form
\Eq{*}{
  y''=\frac{\DET{f''&f\\g''&g}}{\DET{f'&f\\g'&g}}y'-\frac{\DET{f''&f'\\g''&g'}}{\DET{f'&f\\g'&g}}y.
}
After multiplying this equation by $\DET{f'&f\\g'&g}$, and rearranging every term to one side of the
equation, we infer that \eq{DE} is equivalent to
\Eq{EE}{
  \DET{y''&y'&y\\f''&f'&f\\g''&g'&g}=0.
}
It is obvious that the functions $y=f$ and $y=g$ are solutions of \eq{EE} and therefore they also solve \eq{DE}.
\end{proof}

Finally, for a real parameter $p\in\R$, introduce the sine and cosine type functions $S_p,C_p:\R\to\R$ by
\Eq{*}{
  S_p(x):=\begin{cases}
           \sin(\sqrt{-p}x) & \mbox{ if } p<0, \\
           x & \mbox{ if } p=0, \\
           \sinh(\sqrt{p}x) & \mbox{ if } p>0, \\
         \end{cases}\qquad\mbox{and}\qquad
  C_p(x):=\begin{cases}
           \cos(\sqrt{-p}x) & \mbox{ if } p<0, \\
           1 & \mbox{ if } p=0, \\
           \cosh(\sqrt{p}x) & \mbox{ if } p>0. \\
        \end{cases}
}
By basic results on second-order linear homogeneous differential equations, it follows that the functions $S_p$ and $C_p$ constructed above form a fundamental system of solutions for the differential equation
\Eq{*}{
  y''=py.
}

Using these notations, we are now in the position to formulate the main theorem of this paper.

\Thm{MT}{
Let $p\in\R$, let $\varphi:I\to\R_+$ be a positive continuous function and let $(f,g),(h,k)\in\B_0(I)$
such that
\Eq{SC}{
(f,g)\sim(S_p/\varphi,C_p/\varphi) \qquad\mbox{and}\qquad
(h,k)\sim(S_p\cdot\varphi,C_p\cdot\varphi).
}
Then, for all $s,t\in\R_+$, the invariance equation \eq{inv} holds. \\\indent
Conversely, let $(f,g),(h,k)\in\B_4(I)$ and $t,s\in\R_+$ with $t\neq s$ such that the functional
equation \eq{inv} be valid. Then there exist a positive 4 times continuously differentiable function
$\varphi:I\to\R_+$ and a real parameter $p\in\R$ such that the equivalences \eq{SC} are satisfied.
}

\section{Formulas for generalized Wronski determinants}

\Lem{2-}{
\begin{enumerate}[(i)]
\item If $(f,g)\in\B_3(I)$, then
     \Eq{W2}{
\frac{W^{3,0}_{f,g}}{W^{1,0}_{f,g}}=\Phi_{f,g}'+\Phi_{f,g}^2+\Psi_{f,g}\qquad\mbox{and}\qquad
\frac{W^{3,1}_{f,g}}{W^{1,0}_{f,g}}=-\Phi_{f,g}\Psi_{f,g}-\Psi_{f,g}';
}
\item If $(f,g)\in\B_4(I)$, then
\Eq{W3}{
\frac{W^{4,0}_{f,g}}{W^{1,0}_{f,g}}
=\Phi_{f,g}''+3\Phi_{f,g}'\Phi_{f,g}+\Phi_{f,g}^3+2\Phi_{f,g}\Psi_{f,g}+2\Psi_{f,g}'.
}
\end{enumerate}
}
\begin{proof}
Let $(f,g)\in\B_3(I)$. Computing the derivative of $\Phi_{f,g}$ and using $W^{1,1}_{f,g}=0$, we get
\Eq{*}{
\Phi'_{f,g}=\frac{W^{3,0}_{f,g}+W^{2,1}_{f,g}}{W^{1,0}_{f,g}}
 -\frac{W^{2,0}_{f,g}\Big(W^{2,0}_{f,g}+W^{1,1}_{f,g}\Big)}{\Big(W^{1,0}_{f,g}\Big)^2}
  =\frac{W^{3,0}_{f,g}}{W^{1,0}_{f,g}} -\Psi_{f,g}-\Phi^2_{f,g}.
}
Hence the first equality of \eq{W2} follows immediately. In order to show the second equality in \eq{W2},
we differentiate $\Psi_{f,g}$ and use $W^{1,1}_{f,g}=W^{2,2}_{f,g}=0$ to obtain
 \Eq{*}{
 \Psi'_{f,g}=-\frac{W^{3,1}_{f,g}+W^{2,2}_{f,g}}{W^{1,0}_{f,g}}
 +\frac{W^{2,1}_{f,g}\Big(W^{2,0}_{f,g}+W^{1,1}_{f,g}\Big)}{\Big(W^{1,0}_{f,g}\Big)^2}
 =-\frac{W^{3,1}_{f,g}}{W^{1,0}_{f,g}}-\Phi_{f,g}\Psi_{f,g}.
}
Therefore, the second formula in \eq{W2} follows directly.

To prove identity \eq{W3}, let $(f,g)\in\B_4(I)$. Then by differentiating the both sides of
the first equality in \eq{W2}, we get
\Eq{*}{
\frac{W^{4,0}_{f,g}+W^{3,1}_{f,g}}{W^{1,0}_{f,g}}
-\frac{W^{3,0}_{f,g}\Big(W^{2,0}_{f,g}+W^{1,1}_{f,g}\Big)}{\Big(W^{1,0}_{f,g}\Big)^2}
=\Phi_{f,g}''+2\Phi_{f,g}\Phi_{f,g}'+\Psi_{f,g}'.
}
Now, substituting the two formulas of \eq{W2} into this equation, equality \eq{W3} follows after a simple
calculation.
\end{proof}

The following transformation rules will also be needed.

\Lem{2+}{If $(f,g)\in\B_2(I)$ and $\varphi:I\to\R$ is a twice differentiable positive function, then
\Eq{2+}{
  \Phi_{\varphi f,\varphi g}&=\Phi_{f,g}+2\frac{\varphi'}{\varphi},\\
  \Psi_{\varphi f,\varphi g}&=\Psi_{f,g}-\frac{\varphi'}{\varphi} \Phi_{f,g}+\frac{\varphi''}{\varphi}-2\Big(\frac{\varphi'}{\varphi}\Big)^2.
}}

\begin{proof} We prove first that, under the assumptions of the lemma, the following identities are valid:
\Eq{WW}{
  W^{1,0}_{\varphi f,\varphi g}&=\varphi^2 W^{1,0}_{f,g},\\
  W^{2,0}_{\varphi f,\varphi g}&=\varphi^2 W^{2,0}_{f,g}+2\varphi'\varphi W^{1,0}_{f,g},\\
  W^{2,1}_{\varphi f,\varphi g}&=\varphi^2 W^{2,1}_{f,g}+\varphi'\varphi W^{2,0}_{f,g}+\big(2(\varphi')^2-\varphi''\varphi\big) W^{1,0}_{f,g}.
}
Let $F$ denote the vector valued function $\Big(\!\!\begin{array}{c}f\\g\end{array}\!\!\Big)$. Then
\Eq{*}{
  W^{1,0}_{\varphi f,\varphi g}
  =\Det{(\varphi F)'}{\varphi F}
  =\Det{\varphi' F+\varphi F'}{\varphi F}
  =\Det{\varphi' F}{\varphi F}+\Det{\varphi F'}{\varphi F}
  =\varphi^2 W^{1,0}_{f,g}.
}
Similarly,
\Eq{*}{
  W^{2,0}_{\varphi f,\varphi g}
  &=\Det{(\varphi F)''}{\varphi F}
  =\Det{\varphi'' F+2\varphi' F'+\varphi F''}{\varphi F}\\
  &=\Det{\varphi'' F}{\varphi F}+2\Det{\varphi' F'}{\varphi F}+\Det{\varphi F''}{\varphi F}
  =\varphi^2 W^{2,0}_{f,g}+2\varphi'\varphi W^{1,0}_{f,g}.
}
and
\Eq{*}{
  W^{2,1}_{\varphi f,\varphi g}
  =&\Det{(\varphi F)''}{(\varphi F)'}
  =\Det{\varphi'' F+2\varphi' F'+\varphi F''}{\varphi' F+\varphi F'}\\
  =&\Det{\varphi'' F}{\varphi' F}+2\Det{\varphi' F'}{\varphi' F}+\Det{\varphi F''}{\varphi' F}\\
    &+\Det{\varphi'' F}{\varphi F'}+2\Det{\varphi' F'}{\varphi F'}+\Det{\varphi F''}{\varphi F'}\\
  =&\varphi^2 W^{2,1}_{f,g}+\varphi'\varphi W^{2,0}_{f,g}+\big(2(\varphi')^2-\varphi''\varphi\big) W^{1,0}_{f,g}.
}
Dividing the second and third identities by the first one side by side, the formulas stated in \eq{2+} follow immediately.
\end{proof}

\section{Partial derivatives of weighted Bajraktarević means}

To solve the invariance equation \eq{inv}, we need to compute the partial derivatives
\Eq{*}{
\partial_1^iB_{f,g}(x,x;t,s)\qquad\mbox{and} \qquad\partial_1^iB_{h,k}(x,x;s,t)\qquad  (i\in\{1,2,3,4\},
x\in I).
}
According to \lem{sol}, we can apply implicit differentiation for the identity
\Eq{d1}{
  \left|\begin{array}{cc}
  tf(x)+sf(y) & f(B_{f,g}(x,y;t,s)) \\ tg(x)+sg(y) & g(B_{f,g}(x,y;t,s))
  \end{array}\right|=0 \qquad(x,y\in I),
}
which yields the partial derivative of the mean $B_{f,g}(x,y;t,s)$ and, completely analogously,
one can calculate the partial derivatives of $B_{h,k}(x,y;s,t)$.

\Thm{PD}{
Let $t,s\in\R_+$ be fixed. Then the following formulas are valid for all $x\in I$.
\begin{enumerate}[(i)]
 \item If $(f,g)\in\B_1(I)$, then
\Eq{m1}{
 \partial_1B_{f,g}(x,x;t,s)=\frac{t}{t+s};
 }
 \item If $(f,g)\in\B_2(I)$, then
\Eq{m2}{
 \partial_1^2B_{f,g}(x,x;t,s)=\frac{ts}{(t+s)^2}\Phi_{f,g}(x);
 }
 \item If $(f,g)\in\B_3(I)$, then
 \Eq{m3}{
 \partial_1^3 B_{f,g}(x,x;t,s)
 &=\frac{ts}{(t+s)^3}\Big((s-t)(\Phi_{f,g}^2+\Psi_{f,g})+(2t+s)\Phi_{f,g}'\Big)(x);
  }
  \item If $(f,g)\in\B_4(I)$, then
 \Eq{m4}{
 \partial_1^4 B_{f,g}(x,x;t,s)
 &=\frac{ts}{(t+s)^4}\Big((s^2+3ts+3t^2)\Phi_{f,g}''+(3s^2+5ts-5t^2)\Phi_{f,g}'\Phi_{f,g}\\
     &\hspace{1cm}+(s^2-4ts+t^2)(\Phi_{f,g}^3+2\Phi_{f,g}\Psi_{f,g})+2(s^2+ts-t^2)\Psi_{f,g}'\Big)(x).
 }
\end{enumerate}
}
\begin{proof}
For brevity, for $x,y\in I$, denote $M(x,y):=B_{f,g}(x,y;t,s)$ and let $F$
denote the vector valued function $\Big(\!\!\begin{array}{c}f\\g\end{array}\!\!\Big)$. Then, with
these notations, \eq{d1} can be rewritten as
\Eq{d1'}{
  \Det{tF(x)+sF(y)}{F(M(x,y))}=0 \qquad(x,y\in I).
}

If $(f,g)\in\B_1(I)$, then differentiating equation \eq{d1} with respect to the first variable, we get
\Eq{p1}{
  \Det{tF'(x)}{F(M(x,y))}  +\Det{tF(x)+sF(y)}{F'(M(x,y))}\cdot\partial_1 M(x,y)=0.
}
Substituting $x=y$, it follows that
\Eq{*}{
  tW_{f,g}^{1,0}(x)+(t+s)W_{f,g}^{0,1}(x)\cdot\partial_1 M(x,x)=0.
}
Dividing both sides by $(t+s)W_{f,g}^{1,0}(x)\neq0$, this equality yields
\Eq{m1+}{
  \partial_1 B_{f,g}(x,x;t,s)=\partial_1 M(x,x)=\frac{t}{t+s},
}
which is exactly formula \eq{m1}.

Let $(f,g)\in\B_2(I)$, then differentiating equation \eq{p1} with respect to the first variable, we have
\Eq{p2}{
  \Det{tF''(x)}{F(M(x,y))}
  &+2\Det{tF'(x)}{F'(M(x,y))}\partial_1 M(x,y)\\
  &+\Det{tF(x)+sF(y)}{F''(M(x,y))}(\partial_1 M(x,y))^2\\
  &+\Det{tF(x)+sF(y)}{F'(M(x,y))}\partial_1^2 M(x,y)=0.
}
Now putting $x=y$, we arrive at
\Eq{*}{
  tW_{f,g}^{2,0}(x)+(t+s)W_{f,g}^{0,2}(x)(\partial_1 M(x,x))^2
  +(t+s)W_{f,g}^{0,1}(x)\cdot\partial_1^2 M(x,x)=0.
}
Substituting the formula for $\partial_1M(x,x)$ obtained in equation \eq{m1+}, and then dividing the equation side by side by $(t+s)W_{f,g}^{1,0}(x)\neq0$, it follows that
\Eq{*}{
  \frac{t}{t+s}\frac{W_{f,g}^{2,0}(x)}{W_{f,g}^{1,0}(x)}
  +\frac{t^2}{(t+s)^2}\frac{W_{f,g}^{2,0}(x)}{W_{f,g}^{1,0}(x)}
  =\partial_1^2 M(x,x),
}
which yields that
\Eq{m2+}{
  \partial_1^2 B_{f,g}(x,x;t,s)=\partial_1^2 M(x,x)
  =\frac{t(t+s)-t^2}{(t+s)^2}\cdot\frac{W_{f,g}^{2,0}(x)}{W_{f,g}^{1,0}(x)}
  =\frac{ts}{(t+s)^2}\Phi_{f,g}(x).
}

To determine the third-order partial derivative $\partial_1^3B_{f,g}(x,y;t,s)$, assume that $(f,g)\in\B_3(I)$.
Then, by computing the derivative of \eq{p2} with respect to the first variable, we get
\Eq{p3}{
  &\Det{tF'''(x)}{F(M(x,y))}
  +3\Det{tF''(x)}{F'(M(x,y))}\partial_1 M(x,y)\\
  &+3\Det{tF'(x)}{F'(M(x,y))}\partial_1^2 M(x,y)
  +3\Det{tF'(x)}{F''(M(x,y))}\big(\partial_1 M(x,y)\big)^2\\
  &+3\Det{tF(x)+sF(y)}{F''(M(x,y))}\partial_1 M(x,y)\partial_1^2 M(x,y)\\
  &+\Det{tF(x)+sF(y)}{F'''(M(x,y))}\big(\partial_1 M(x,y)\big)^3
  +\Det{tF(x)+sF(y)}{F'(M(x,y))}\partial_1^3 M(x,y)=0.
}
Taking $x=y$ and using formulae \eq{m1+} and \eq{m2+}, hence
\Eq{*}{
  t&W_{f,g}^{3,0}(x)+\frac{3t^2}{t+s}W_{f,g}^{2,1}(x)
  +\frac{3t^3}{(t+s)^2}W_{f,g}^{1,2}(x)+\frac{3t^2s}{(t+s)^2}W_{f,g}^{0,2}(x)\Phi_{f,g}(x)
  +\frac{t^3}{(t+s)^2}W_{f,g}^{0,3}(x)\\
  &+(t+s)W_{f,g}^{0,1}(x)\cdot\partial_1^3 M(x,x)=0
}
After simple calculations, the last equation reduces to
\Eq{*}{
  \partial_1^3 M(x,x)
  =\frac{ts(2t+s)}{(t+s)^3}\cdot\frac{W_{f,g}^{3,0}(x)}{W_{f,g}^{1,0}(x)}
  +\frac{3t^2s}{(t+s)^3}\cdot\frac{W_{f,g}^{2,1}(x)}{W_{f,g}^{1,0}(x)}
  -\frac{3t^2s}{(t+s)^3}\frac{W_{f,g}^{2,0}(x)}{W_{f,g}^{1,0}(x)}\Phi_{f,g}(x)
}
Using the first equation in \eq{W2} and the the definitions of $\Phi_{f,g}$ and $\Psi_{f,g}$, the following formula follows
 \Eq{m3+}{
 \partial_1^3 B_{f,g}(x,x;t,s)=\partial_1^3 M(x,x)
 =\frac{ts}{(t+s)^3}\Big((s-t)(\Phi_{f,g}^2+\Psi_{f,g})+(2t+s)\Phi_{f,g}'\Big)(x).
}

Finally, assume that $(f,g)\in\B_4(I)$. Differentiating equation \eq{p3} with respect to $x$, we get
\Eq{*}{
&\Det{tF''''(x)}{F(M(x,y))}
  +4\Det{tF'''(x)}{F'(M(x,y))}\partial_1 M(x,y)\\
  &+6\Det{tF''(x)}{F''(M(x,y))}\big(\partial_1 M(x,y)\big)^2
  +6\Det{tF''(x)}{F'(M(x,y))}\partial_1^2 M(x,y)\\
  &+4\Det{tF'(x)}{F'(M(x,y))}\partial_1^3 M(x,y)
  +12\Det{tF'(x)}{F''(M(x,y))}\partial_1 M(x,y)\partial_1^2 M(x,y)\\
  &+4\Det{tF'(x)}{F'''(M(x,y))}\big(\partial_1 M(x,y)\big)^3
  +\Det{tF(x)+sF(y)}{F''''(M(x,y))}\big(\partial_1 M(x,y)\big)^4\\
  &+6\Det{tF(x)+sF(y)}{F'''(M(x,y))}\big(\partial_1 M(x,y)\big)^2\partial_1^2 M(x,y)\\
  &+\Det{tF(x)+sF(y)}{F''(M(x,y))}\Big(3\big(\partial_1^2 M(x,y)\big)^2+4\partial_1 M(x,y)\partial_1^3 M(x,y)\Big)\\
  &+\Det{tF(x)+sF(y)}{F'(M(x,y))}\partial_1^4 M(x,y)=0.
}
Analogously, putting $x=y$, using \eq{m1+}, \eq{m2+}, and \eq{m3+}, we obtain the following equality
\Eq{*}{
  t&W_{f,g}^{4,0}(x)
  +\frac{4t^2}{t+s}W_{f,g}^{3,1}(x)
  +\frac{6t^2s}{(t+s)^2}W_{f,g}^{2,1}(x)\Phi_{f,g}(x)
  +\frac{12t^3s}{(t+s)^3}W_{f,g}^{1,2}(x)\Phi_{f,g}(x)\\
  &+\frac{4t^4}{(t+s)^3}W_{f,g}^{1,3}(x)
  +\frac{t^4}{(t+s)^3}W_{f,g}^{0,4}(x)
  +\frac{6t^3s}{(t+s)^3}W_{f,g}^{0,3}(x)\Phi_{f,g}(x)
  +\frac{3t^2s^2}{(t+s)^3}W_{f,g}^{0,2}(x)\Phi_{f,g}^2(x)\\
  &+\frac{4t^2s}{(t+s)^3}W_{f,g}^{0,2}(x)\Big((s-t)(\Phi_{f,g}^2+\Psi_{f,g})+(2t+s)\Phi_{f,g}'\Big)(x)
  +(t+s)W_{f,g}^{0,1}(x)\cdot\partial_1^4 M(x,x)=0.
}
Therefore,
\Eq{*}{
  \partial_1^4 M(x,x)=&
  \frac{ts}{(t+s)^4}\Bigg((3t^2+3st+s^2)\frac{W_{f,g}^{4,0}}{W_{f,g}^{1,0}}
  +4t(2t+s)\frac{W_{f,g}^{3,1}}{W_{f,g}^{1,0}}
  -6t^2\frac{W_{f,g}^{3,0}}{W_{f,g}^{1,0}}\Phi_{f,g}\\
  &+6t(s-t)\frac{W_{f,g}^{2,1}}{W_{f,g}^{1,0}}\Phi_{f,g}
  -t\cdot\frac{W_{f,g}^{2,0}}{W_{f,g}^{1,0}}\Big((7s-4t)\Phi_{f,g}^2+4(s-t)\Psi_{f,g}+4(2t+s)\Phi_{f,g}'\Big)\Bigg)(x).
}
Using the identities in \eq{W2} and \eq{W3}, the desired formula \eq{m4} follows immediately.
\end{proof}

\section{Necessity}
\Lem{N1}{
Let $t,s\in\R_+$ and let $(f,g),(h,k)\in\B_2(I)$ satisfy the functional equation \eq{inv}. Then
\Eq{C1}{
\Phi_{f,g}+\Phi_{h,k}=0
}
\begin{proof}
Differentiation equation \eq{inv} twice with respect to $x$, then substituting $y=x$, we get
 \Eq{*}{
 \partial_1^2 B_{f,g}(x,x;t,s)+\partial_1^2 B_{h,k}(x,x;s,t)=0 \qquad(x\in I).
}
Then by the second-order partial derivative formula in \eq{m2}, we obtain
\Eq{*}{
\frac{ts}{(t+s)^2}\Phi_{f,g}(x)+\frac{ts}{(t+s)^2}\Phi_{h,k}(x)=0 \qquad(x\in I),
}
which implies the equality \eq{C1}.
\end{proof}
}
Assuming that $(f,g)\in\B_2(I)$, for simplicity, denote $\Phi_{f,g}$ by $\Phi$. Then, if \eq{inv} holds, according to the previous lemma, we have that
\Eq{Phi}{
\Phi_{h,k}=-\Phi.
}

\Lem{N2}{
Let $t,s\in\R_+$ with $t\neq s$ and let $(f,g),(h,k)\in\B_3(I)$ satisfy the functional equation \eq{inv}. Then
\Eq{C2}{
\Phi'+\Psi_{h,k}-\Psi_{f,g}=0.
}
\begin{proof}
Differentiating \eq{inv} three times with respect to the variable $x$, then substituting $y=x$, we get
\Eq{*}{
\partial_1^3 B_{f,g}(x,x;t,s)+\partial_1^3 B_{h,k}(x,x;s,t)=0 \qquad(x\in I).
}
Then the formula for the third-order partial derivatives by \eq{m3} implies that
\Eq{*}{
\big((s-t)(\Phi_{f,g}^2+\Psi_{f,g})+(2t+s)\Phi_{f,g}'\big)
 +\big((t-s)(\Phi_{h,k}^2+\Psi_{h,k})+(2s+t)\Phi_{h,k}'\big)=0.
}
Using $t\neq s$ and \eq{Phi}, the last equation reduces to the desired equality \eq{C2}.
\end{proof}
}

\Lem{N3}{
Let $t,s\in\R_+$ with $t\neq s$ and let $(f,g),(h,k)\in\B_4(I)$ satisfy the functional equation \eq{inv}.
Then there exists a real constant $p$ such that
\Eq{C3}{
\Psi_{f,g}=\frac12\Phi'-\frac14\Phi^2+p \qquad\mbox{and}\qquad
  \Psi_{h,k}=-\frac12\Phi'-\frac14\Phi^2+p.
}}

\begin{proof}
Differentiating \eq{inv} four times with respect to the variable $x$, then substituting $y=x$, we arrive at
\Eq{*}{
\partial_1^4 B_{f,g}(x,x;t,s)+\partial_1^4 B_{h,k}(x,x;s,t)=0 \qquad(x\in I).
}
Then the formula for the fourth-order partial derivatives by \eq{m3} and equality \eq{Phi} yield that
\Eq{*}{
(t^2-s^2)\Phi''&-(t^2-5ts+s^2)\Phi'\Phi+(s^2-4ts+t^2)\Phi(\Psi_{f,g}-\Psi_{h,k})\\
                    &+(s^2+ts-t^2)\Psi_{f,g}'+(t^2+ts-s^2)\Psi_{h,k}'=0.
}
After simple calculations and using formula \eq{C2}, we get
\Eq{*}{
  (\Phi^2)'+2(\Psi'_{h,k}+\Psi'_{f,g})=0.
 }
By integrating the last equation, we get that there exists a real constant $p$ such that
\Eq{*}{
  \Phi^2+2(\Psi_{h,k}+\Psi_{f,g})=4p.
 }
This equality, combined with \eq{C2}, implies the formulas stated in \eq{C3}.
\end{proof}

Finally, we can complete the proof of \thm{MT}.

\begin{proof}[Proof of the necessity in \thm{MT}]
Let $(f,g),(h,k)\in\B_4(I)$ and $t,s\in\R_+$ with $t\neq s$ such that the functional equation \eq{inv} be valid.
Then, as we have seen it in \lem{N3}, there exist a real constant $p\in\R$ such that \eq{C3} holds. Define
\Eq{*}{
  \varphi:=\frac1{\sqrt{|W^{1,0}_{f,g}|}}=|W^{1,0}_{f,g}|^{-\frac12}.
}
Then
\Eq{*}{
  \frac{\varphi'}{\varphi}=-\frac{(W^{1,0}_{f,g})'}{2W^{1,0}_{f,g}}
  =-\frac{W^{2,0}_{f,g}}{2W^{1,0}_{f,g}}=-\frac12\Phi_{f,g}=-\frac12\Phi.
}
Differentiating again, it follows that
\Eq{*}{
  \frac{\varphi''}{\varphi}=-\frac12\Phi'+\frac14\Phi^2.
}
Now, using \lem{2+} and the first formula from \eq{C3}, we get
\Eq{*}{
  \Phi_{\varphi f,\varphi g}=\Phi_{f,g}+2\frac{\varphi'}{\varphi}=\Phi+2\Big(-\frac12\Phi\Big)=0
}
and
\Eq{*}{
  \Psi_{\varphi f,\varphi g}
  &=\Psi_{f,g}-\frac{\varphi'}{\varphi}\Phi_{f,g}+\frac{\varphi''}{\varphi}-2\Big(\frac{\varphi'}{\varphi}\Big)^2 \\
  &=\Big(\frac12\Phi'-\frac14\Phi^2+p\Big)+\frac12\Phi^2-\frac12\Phi'+\frac14\Phi^2-2\Big(-\frac12\Phi\Big)^2=p.
}
In view of \lem{DE}, these identities imply that $\varphi f$ and $\varphi g$ are solutions of the second-order homogeneous linear differential equation $y''=py$. Hence, $(\varphi f,\varphi g)\sim(S_p,C_p)$, which implies $(f,g)\sim(S_p/\varphi,C_p/\varphi)$.

A completely analogous argument shows that $(h,k)\sim(S_p\cdot\varphi,C_p\cdot\varphi)$. Thus, we have proved that all the solutions of \eq{inv} are of the form stated in \thm{MT}.

\end{proof}

\section{Sufficiency}

\begin{proof}[Proof of the sufficiency in \thm{MT}]
Assume that $p\in\R$ and $\varphi:I\to\R_+$ is a continuous function and \eq{SC} holds on $I$.

If $p\geq0$ then $C_p$ is a positive function. Thus, by \lem{equ0}, we have that $B_{f,g}=B_{S_p/\varphi,C_p/\varphi}$ and $B_{h,k}=B_{S_p\cdot\varphi,C_p\cdot\varphi}$, therefore, it suffices to show that
\Eq{NIE}{
  B_{S_p/\varphi,C_p/\varphi}(x,y;t,s)+B_{S_p\cdot\varphi,C_p\cdot\varphi}(x,y;s,t)=x+y \qquad(x,y\in I).
}

In the case $p=0$, this equation can be rewritten as
\Eq{*}{
  \frac{t\frac{x}{\varphi(x)}+s\frac{y}{\varphi(y)}}{t\frac{1}{\varphi(x)}+s\frac{1}{\varphi(y)}}
  +\frac{sx\varphi(x)+ty\varphi(y)}{s\varphi(x)+t\varphi(y)}=x+y,
}
which is an easy-to-see identity.

For the case $p>0$, we recall the addition theorem for the inverse of the tangent hyperbolic function, which reads as follows:
\Eq{tanh}{
  \tanh^{-1}(a)+ \tanh^{-1}(b)= \tanh^{-1}\Big(\frac{a+b}{1+ab}\Big) \qquad(a,b\in]-1,1[).
}
We can easily see that
\Eq{*}{
  B_{S_p/\varphi,C_p/\varphi}(x,y;t,s)
  &=\frac1{\sqrt{p}}\tanh^{-1}\bigg(\frac{t\varphi(y)\sinh(\sqrt{p}x)+s\varphi(x)\sinh(\sqrt{p}y)}
  {t\varphi(y)\cosh(\sqrt{p}x)+s\varphi(x)\cosh(\sqrt{p}y)}\bigg),\\
  B_{S_p\cdot\varphi,C_p\cdot\varphi}(x,y;s,t)
  &=\frac1{\sqrt{p}}\tanh^{-1}\bigg(\frac{s\varphi(x)\sinh(\sqrt{p}x)+t\varphi(y)\sinh(\sqrt{p}y)}
  {s\varphi(x)\cosh(\sqrt{p}x)+t\varphi(y)\cosh(\sqrt{p}y)}\bigg).
}
Therefore, with the substitutions $u:=\sqrt{p}x$, $v:=\sqrt{p}y$, and $T:=t\varphi(y)$, $S:=s\varphi(x)$, and using the addition formula \eq{tanh}, the left hand side of equation \eq{NIE} can be rewritten and calculated as
\Eq{*}{
  \frac1{\sqrt{p}}&\bigg(\tanh^{-1}\Big(\frac{T\sinh(u)+S\sinh(v)}{T\cosh(u)+S\cosh(v)}\Big)
  +\tanh^{-1}\Big(\frac{S\sinh(u)+T\sinh(v)}{S\cosh(u)+T\cosh(v)}\Big)\bigg)\\
  &=\frac1{\sqrt{p}}\tanh^{-1}\left(\frac{\frac{T\sinh(u)+S\sinh(v)}{T\cosh(u)+S\cosh(v)}+\frac{S\sinh(u)+T\sinh(v)}{S\cosh(u)+T\cosh(v)}}
    {1+\frac{T\sinh(u)+S\sinh(v)}{T\cosh(u)+S\cosh(v)}\frac{S\sinh(u)+T\sinh(v)}{S\cosh(u)+T\cosh(v)}}\right)\\
  &=\frac1{\sqrt{p}}\tanh^{-1}\left(\frac{(T^2+S^2)\sinh(u+v)+TS(\sinh(2u)+\sinh(2v))}
    {(T^2+S^2)\cosh(u+v)+TS(\cosh(2u)+\cosh(2v))}\right)\\
  &=\frac1{\sqrt{p}}\tanh^{-1}\left(\frac{(T^2+S^2)\sinh(u+v)+2TS\sinh(u+v)\cosh(u-v)}
    {(T^2+S^2)\cosh(u+v)+2TS\cosh(u+v)\cosh(u-v)}\right)\\
  &=\frac1{\sqrt{p}}\tanh^{-1}\big(\tanh(u+v)\big)\\
  &=\frac{u+v}{\sqrt{p}}=x+y.
}
This completes the proof of the identity \eq{NIE}.

For the discussion of the case $p<0$, we recall the folk-addition theorem for the inverse of the tangent function, which reads as follows:
\Eq{tan}{
  \tan^{-1}(a)+ \tan^{-1}(b)
  = \begin{cases}
     -\pi+\tan^{-1}\big(\frac{a+b}{1-ab}\big) &\mbox{if } ab>1,\,a,b<0,\\[2mm]
     -\tfrac{\pi}{2} &\mbox{if } ab=1,\,a,b<0, \\[2mm]
     \tan^{-1}\big(\frac{a+b}{1-ab}\big) &\mbox{if } ab<1,\\[2mm]
     \tfrac{\pi}{2} &\mbox{if } ab=1,\,a,b>0, \\[2mm]
     \pi+\tan^{-1}\big(\frac{a+b}{1-ab}\big) &\mbox{if } ab>1,\,a,b>0.
    \end{cases}
}

In the case $p<0$, the functions $C_p$ is not necessarily nonvanishing on $I$, therefore, we cannot immediately express the means $B_{f,g}$ and $B_{h,k}$ in terms of $C_p$ and $S_p$. For the sake of brevity, denote $\sqrt{-p}$ by $q$. First, we shall show that there exist $\alpha\in]-\pi,\pi]$ and $k\in\Z$ such that
\Eq{incl}{
  I\subseteq \big]\tfrac1q\big(\alpha+\big(2k-\tfrac12\big)\pi\big),\tfrac1q\big(\alpha+\big(2k+\tfrac12\big)\pi\big)\big[.
}
In view of the equivalence $(f,g)\sim(S_p/\varphi,C_p/\varphi)$ we have that $(f\varphi,g\varphi)\sim(S_p,C_p)$, hence there exist real constants $c,d$ such that $g\varphi=cS_p+dC_p$ holds on $I$. The function $g$ is nowhere zero on $I$, therefore $cS_p+dC_p$ is also nonvanishing on $I$. We may assume that this function is positive on $I$ (otherwise we replace $(c,d)$ by $(-c,-d)$).
Determine the unique $\alpha\in]-\pi,\pi]$ such that the equalities
\Eq{*}{
  \sin\alpha=\frac{c}{\sqrt{c^2+d^2}},\qquad \cos\alpha=\frac{d}{\sqrt{c^2+d^2}}
}
be valid. Then, we have that $\sin(\alpha) S_p+\cos(\alpha) C_p$ is positive on $I$, that is, for all $x\in I$,
\Eq{*}{
  0<\sin(\alpha) S_p(x)+\cos(\alpha) C_p(x)=\cos(qx-\alpha).
}
With the notation $J:=qI-\alpha$ and the substitution $w:=qx-\alpha$, this condition means that, for all $w\in J$, $\cos(w)>0$. Because $J$ is an interval, it follows that there exists $k\in\Z$ such that
\Eq{*}{
J\subseteq\big]\big(2k-\tfrac12\big)\pi,\big(2k+\tfrac12\big)\pi\big[.
}
This inclusion directly implies \eq{incl}.

For the computation of the means $B_{f,g}$ and $B_{h,k}$, observe first that
\Eq{*}{
(S_p,C_p)
   &\sim \big(\cos(\alpha) S_p-\sin(\alpha) C_p,\sin(\alpha) S_p+\cos(\alpha) C_p\big)\\
   &=\big(\sin(q(\cdot)-\alpha-2k\pi),\cos(q(\cdot)-\alpha-2k\pi)\big).
}
Therefore, by the equivalences, $(f\varphi,g\varphi)\sim(S_p,C_p)$ and $(h/\varphi,k/\varphi)\sim(S_p,C_p)$, it follows that
\Eq{*}{
(f,g)&\sim\big(\tfrac1\varphi\sin(q(\cdot)-\alpha-2k\pi),\tfrac1\varphi\cos(q(\cdot)-\alpha-2k\pi)\big),\\
(h,k)&\sim\big(\varphi\sin(q(\cdot)-\alpha-2k\pi),\varphi\cos(q(\cdot)-\alpha-2k\pi)).
}
Since the function $\cos(q(\cdot)-\alpha-2k\pi)$ is positive on $I$, using \lem{equ0}, we have that
\Eq{*}{
  B_{f,g}(x,y;t,s)
  &=\frac1{q}\bigg(\tan^{-1}\bigg(\frac{\tfrac{t}{\varphi(x)}\sin(qx-\alpha-2k\pi)+\tfrac{s}{\varphi(y)}\sin(qy-\alpha-2k\pi)}
  {\tfrac{t}{\varphi(x)}\cos(qx-\alpha-2k\pi)+\tfrac{s}{\varphi(y)}\cos(qy-\alpha-2k\pi)}\bigg)+\alpha+2k\pi\bigg),\\
  B_{h,k}(x,y;s,t)
  &=\frac1{q}\bigg(\tan^{-1}\Big(\frac{s\varphi(x)\sin(qx-\alpha-2k\pi)+t\varphi(y)\sin(qy-\alpha-2k\pi)}
  {s\varphi(x)\cos(qx-\alpha-2k\pi)+t\varphi(y)\cos(qy-\alpha-2k\pi)}\Big)+\alpha+2k\pi\bigg).
}
Now, denote $u:=qx-\alpha-2k\pi$, $v:=qy-\alpha-2k\pi$, and $T:=t\varphi(y)$, $S:=s\varphi(x)$. Then, by the inclusion \eq{incl}, we have
$u,v\in\big]-\frac\pi2,\frac\pi2\big[$. Hence, it follows that $u+v,u-v\in]-\pi,\pi[$.

With these substitutions, the left hand side of equation \eq{NIE} can be rewritten as
\Eq{*}{
  L:=\frac1{q}\bigg(\tan^{-1}\Big(\frac{T\sin(u)+S\sin(v)}{T\cos(u)+S\cos(v)}\Big)+\alpha+2k\pi
  +\tan^{-1}\Big(\frac{S\sin(u)+T\sin(v)}{S\cos(u)+T\cos(v)}\Big)+\alpha+2k\pi\bigg).
}
In order to calculate $L$ and to show that $L=x+y$, we have to distinguish several cases.

First we consider the case when $T=S$.
Now, observing that $\tfrac{u+v}2,\tfrac{u-v}2\in\big]-\frac\pi2,\frac\pi2\big[$, we have $\cos\big(\tfrac{u-v}2\big)>0$ and
$\tan^{-1}\big(\tan\big(\tfrac{u+v}2\big)\big)=\tfrac{u+v}2$, therefore,
\Eq{*}{
  L&=\frac2{q}\tan^{-1}\Big(\frac{\sin(u)+\sin(v)}{\cos(u)+\cos(v)}\Big)+\frac{2\alpha+4k\pi}{q}
    =\frac2{q}\tan^{-1}\bigg(\frac{2\sin\big(\frac{u+v}2\big)\cos\big(\frac{u-v}2\big)}
   {2\cos\big(\frac{u+v}2\big)\cos\big(\frac{u-v}2\big)}\bigg)+\frac{2\alpha+4k\pi}{q}\\
   &=\frac2{q}\tan^{-1}\bigg(\frac{\sin\big(\frac{u+v}2\big)}
   {\cos\big(\frac{u+v}2\big)}\bigg)+\frac{2\alpha+4k\pi}{q}
   =\frac2{q}\cdot\frac{u+v}{2}+\frac{2\alpha+4k\pi}{q}=\frac{u+v+2\alpha+4k\pi}{q}=x+y.\\
}
From now on, we may assume that $T\neq S$. In this case,
\Eq{pos}{
  T^2+S^2+2TS\cos(u-v)=(T-S)^2+2TS(1+\cos(u-v))>0.
}
Let us denote
\Eq{*}{
  a:=a(u,v):=\frac{T\sin(u)+S\sin(v)}{T\cos(u)+S\cos(v)},\qquad b:=b(u,v):=\frac{S\sin(u)+T\sin(v)}{S\cos(u)+T\cos(v)}.
}
With this notation, we can rewrite $L$ in the following form
\Eq{*}{
  L:=\frac1{q}\big(\tan^{-1}(a)+\tan^{-1}(b)+2\alpha+4k\pi\big).
}
In order to apply the addition formula \eq{tan}, we have to clarify to which subdomain the pair $(a,b)$ belongs to. We shall prove that
\Eq{*}{
  u+v &\in\big]-\pi,-\tfrac{\pi}{2}\big[\qquad &\Longrightarrow& \qquad ab>1,\,a,b<0,\\
  u+v &=-\tfrac{\pi}{2}\qquad &\Longrightarrow& \qquad ab=1,\,a,b<0,\\
  u+v &\in\big]-\tfrac{\pi}{2},\tfrac{\pi}{2}\big[\qquad &\Longrightarrow& \qquad ab<1,\\
  u+v &=\tfrac{\pi}{2}\qquad &\Longrightarrow& \qquad ab=1,\,a,b>0,\\
  u+v &\in\big]\tfrac{\pi}{2},\pi\big[\qquad &\Longrightarrow& \qquad ab>1,\,a,b>0.
}
By $u,v\in \big]-\tfrac{\pi}{2},\tfrac{\pi}{2}\big[$, we always have $\cos(u),\cos(v)>0$.
If $u+v \in\big]-\pi,-\tfrac{\pi}{2}\big]$, then $u,v\in \big]-\tfrac{\pi}{2},0\big[$, which implies that $\sin(u),\sin(v)<0$. Therefore, $a,b<0$. Similarly, if $u+v \in\big[\tfrac{\pi}{2},\pi\big[$, then $a,b>0$ holds. A direct computation shows that
\Eq{*}{
  1-ab=\frac{(T^2+S^2+2TS\cos(u-v))\cos(u+v)}{(T\cos(u)+S\cos(v))(S\cos(u)+T\cos(v))}.
}
Therefore, using \eq{pos}, the sign of $1-ab$ equals the sign of $\cos(u+v)$. This implies $ab<1$, $ab=1$, and $ab>1$ if and only if
$u+v\in\big]-\tfrac{\pi}{2},\tfrac{\pi}{2}\big[$, $u+v=\pm\tfrac\pi2$, and $u+v\in ]-\pi,-\tfrac{\pi}{2}\big[\cup \big]\tfrac{\pi}{2},\pi\big[$, respectively.

In order to complete the proof of sufficiency, we shall calculate $L$ under consideration all cases stated  in \eq{tan}.
For brevity, define
\Eq{*}{
  \sigma:=\sigma(u,v):=
    \begin{cases}
     -1 &\mbox{if } u+v \in\big]-\pi,-\tfrac{\pi}{2}\big],\\[2mm]
     0&\mbox{if } u+v \in\big]-\tfrac{\pi}{2},\tfrac{\pi}{2}\big[,\\[2mm]
     1 &\mbox{if } u+v \in\big[\tfrac{\pi}{2},\pi\big[.
    \end{cases}
}
In the case $u+v\neq \pm\frac\pi2$ we have that $ab\neq1$. Therefore, using the addition formula \eq{tan}, we get
\Eq{*}{
  L&=\frac1{q}\tan^{-1}\left(\frac{a+b}{1-ab}\right)+\frac{2\alpha+4k\pi+\sigma\pi}{q}\\
   &=\frac1{q}\tan^{-1}\left(\frac{\frac{T\sin(u)+S\sin(v)}{T\cos(u)+S\cos(v)}+\frac{S\sin(u)+T\sin(v)}{S\cos(u)+T\cos(v)}}
    {1-\frac{T\sin(u)+S\sin(v)}{T\cos(u)+S\cos(v)}\frac{S\sin(u)+T\sin(v)}{S\cos(u)+T\cos(v)}}\right)+\frac{2\alpha+4k\pi+\sigma\pi}{q}\\
  &=\frac1{q}\tan^{-1}\left(\frac{(T^2+S^2)\sin(u+v)+TS(\sin(2u)+\sin(2v))}
    {(T^2+S^2)\cos(u+v)+TS(\cos(2u)+\cos(2v))}\right)+\frac{2\alpha+4k\pi+\sigma\pi}{q}\\
  &=\frac1{q}\tan^{-1}\left(\frac{(T^2+S^2)\sin(u+v)+2TS\sin(u+v)\cos(u-v)}
    {(T^2+S^2)\cos(u+v)+2TS\cos(u+v)\cos(u-v)}\right)+\frac{2\alpha+4k\pi+\sigma\pi}{q}\\
  &=\frac1{q}\tan^{-1}\big(\tan(u+v)\big)+\frac{2\alpha+4k\pi+\sigma\pi}{q}=\frac{u+v-\sigma\pi}{q}+\frac{2\alpha+4k\pi+\sigma\pi}{q}\\
  &=\frac{u+v+2\alpha+4k\pi}{q}=x+y.
}
Finally assume that $u+v=\sigma\frac\pi2$. Then $ab=1$ such that either $a,b<0$ or $a,b>0$. Therefore, by the addition formula \eq{tan} again, it results that
\Eq{*}{
\tan^{-1}(a)+ \tan^{-1}(b)=\sigma\frac\pi2=u+v.
}
Using this equality, we can see at once that $L$ is equal to $x+y$ also in this case. This completes the proof of the sufficiency in \thm{MT}.
\end{proof}

\def\MR{}


\end{document}